\documentclass{article}
\usepackage{amsmath,amssymb}

\begin{document}

\title{Characterization of Absolutely Norm Attaining Compact Hyponormal Operators\thanks{%
Mathematics Subject Classifications:  47A75, 47A10}}
\date{{\small Received 20 October 2018}}
\author{Benard Okelo\thanks{%
Institute  of Mathematics, University of Muenster, Einstein Street 62, 48149-Muenster, Germany.}}
\maketitle

\begin{abstract}
In this paper, we characterize absolute norm-attainability for compact hyponormal operators. We give necessary and sufficient conditions for a bounded linear compact hyponormal operator  on an infinite dimensional complex Hilbert space to be absolutely norm attaining. Moreover, we discuss the structure of compact hyponormal operators when they are self-adjoint and normal. Lastly,  we discuss in general, other properties of compact hyponormal operators when they are absolutely norm attaining and their commutators.
\end{abstract}

\section{Introduction}

The study of norm attaining operators has been interesting to many mathematicians and researchers over decades( see \cite{Car}, \cite{Enf} and \cite{Ram}). The class of absolutely norm attaining operators between complex Hilbert spaces was introduced by  \cite{Car} and they discussed several important  examples and properties of these operators. The class of absolutely norm attaining operators is denoted by $\mathcal{AN}(H).$ A synonymous class called norm-attainable operators have also been discussed by Okelo in \cite{Oke} and it has been determined that they share similar characteristics. In this paper,  we give  necessary and sufficient conditions for an operator to be hyponormal  and belongs to $\mathcal{AN}(H)$. In fact, we show that  a bounded operator $T$ defined on an infinite dimensional Hilbert space is hyponormal and belongs to $\mathcal{AN}(H)$ if and only if there exists a unique triple $(K,F,\alpha)$, where $K$ is a positive compact operator, $F$ is a positive finite rank operator, $\alpha$ is positive real number such that  $T=K-F+\alpha I$ and $KF=0,\; F\leq \alpha I$. In fact, here $\alpha=m_e(T)$, the essential minimum modulus of $T$.  Moreover,  we give explicit structure of self-adjoint  and $\mathcal{AN}$-operators as well as   hyponormal and  $\mathcal{AN}$-operators. Finally, we also obtain structure of general $\mathcal{AN}$-operators. In the process we also prove several important properties of $\mathcal{AN}$-operators. Unless otherwise stated, the hyponormal operators in this work are compact. We organize the article as follows: Section 1: Introduction; Section 2:  Preliminaries and notations; Section 3: Main results.
\section{Preliminaries}
In this section, we give the preliminaries. These include the basic terms, definitions and notations which are useful in the sequel. Throughout the paper, we consider  all Hilbert spaces to be infinite dimensional and complex.   We denote inner product and the
induced norm by  by  $\langle . , .\rangle$ and $\|\cdot \|$ respectively. The unit sphere  of a closed subspace $M$ of $H$ is denoted by $S_M:={\{x\in M: \|x\|=1}\}$ and  $P_M$ denote the orthogonal projection $P_M:H\rightarrow H$ with range $M$. The identity operator on $M$ is denoted by $I_M$. See details in \cite{Car} and the references therein.\\

DEFINITION 2.1. An operator $T: H_1\rightarrow H_2$ is said to be bounded if there exists a $C>0$ such that $\|Tx\|\leq C \|x\|,$ for all $x\in H_1$.  If $T$ is bounded, the quantity $\|T\|=\sup{\{\|Tx\|: x\in S_{H_1}}\}$ is finite and is called the norm of $T$.\\

 We denote the space of all bounded linear operators between $H_1$ and $H_2$  by $\mathcal B(H_1,H_2)$. In general, the set of all bounded linear operators on $H$ is denoted by $\mathcal B(H)$.\\

DEFINITION 2.2. For $T\in \mathcal B(H_1,H_2)$, there exists a unique operator denoted by $T^*:H_2\rightarrow H_1$ called the adjoint operator satisfying $\langle Tx,y\rangle =\langle x, T^*y\rangle, \; \text{for all}\; x\in H_1 \; \text{and} \; \text{for all}\; y\in H_2.$\\

DEFINITION 2.3. An operator $T\in \mathcal B(H_1,H_2)$ is said to be norm attaining if there exists a $x\in S_{H_1}$ such that $\|Tx\|=\|T\|$.  We denote the class of norm attaining operators by $\mathcal N(H_1,H_2)$. \\

REMARK 2.1. It is known that $\mathcal N(H_1,H_2)$ is dense in $\mathcal B(H_1,H_2)$ with respect to the operator norm of $\mathcal B(H_1,H_2)$. We refer to \cite{Enf} for more details.\\

DEFINITION 2.4(\cite{Car}). An operator $T\in \mathcal B(H_1,H_2)$ is said to be  absolutely norm attaining  or $\mathcal {AN}$-operator (shortly), if  $T|_M$,  the restriction of $T$ to $M $,  is norm attaining for every non zero closed subspace $M$ of $H_1$. That is $T|_M\in \mathcal N(M, H_2)$ for every non zero closed subspace $M$ of $H_1.$  \\

DEFINITION 2.5. An operator $T\in \mathcal{B}(H)$ is said to be hyponormal if $\|T^{*}x\|\leq\|Tx\|, $ for all $x\in H.$\\

REMARK 2.2. This class contains $\mathcal K(H_1,H_2)$, and  the class of partial isometries with finite dimensional null space or finite dimensional range space. \\

In the remaining part of this section, we give standard terminologies and notations found in \cite{Hal}. Let $T\in \mathcal B(H)$. Then $T$ is said to be normal if $T^*T=TT^*$, self-adjoint if $T=T^*$. If $T$ is self-adjoint and $\langle Tx,x\rangle \geq 0,$ for all $x\in H$, then $T$ is called positive. It is well known that for a positive operator $T$, there exists a unique  positive operator $S\in \mathcal B(H)$ such that $S^2=T$. We write $S=T^{\frac{1}{2}}$  and is called as the positive square root of $T$. If $S,T\in \mathcal B(H)$ are self-adjoint and $S-T\geq 0$, then we write this by $S\geq T$. If $P\in \mathcal B(H)$ is such that $P^2=P$, then $P$ is called a projection. If Null space of $P$,  $N(P)$ and range of $P,$ $R(P)$ are orthogonal to each other, then $P$ is called an orthogonal projection. It is a well known  fact  that a projection  $P$ is an orthogonal projection if and only if it is self-adjoint if and only if it is normal. We call an operator $V\in \mathcal B(H_1,H_2)$ to be an isometry if $\|Vx\|=\|x\|,$ for each $x\in H_1$. An operator $V\in \mathcal B(H_1,H_2)$ is said to be a partial isometry if $V|_{N(V)^{\bot}}$ is an isometry. That is, $\|Vx\|=\|x\|$ for all $x\in N(V)^{\bot}$. If $V\in \mathcal B(H)$ is isometry and onto, then $V$ is said to be a unitary operator. If $T\in \mathcal B(H) $ is a self-adjoint operator, then $T=T_{+}-T_{-}$, where $T_{\pm}$ are positive operators. Here $T_{+}$ is called the positive part and $T_{-}$ is called the negative part of $T$ . Moreover, this decomposition is unique. In general, if $T\in \mathcal B(H_1,H_2)$, then $T^*T\in \mathcal B(H_1)$ is positive and $|T|:=(T^*T)^{\frac{1}{2}}$ is called the modulus of $T$. In fact, there exists a unique partial isometry $V\in \mathcal B(H_1,H_2)$ such that $T=V|T|$ and $N(V)=N(T)$. This factorization is called the polar decomposition of $T$. If $T\in \mathcal B(H)$, then $T=\frac{T+T^*}{2}+i(\frac{T-T^*}{2i})$. The operators $Re(T):=\frac{T+T^*}{2}$ and $Im(T):=\frac{T-T^*}{2i}$ are self-adjoint and called the real and the imaginary parts of $T$ respectively. A closed subspace $M$ of $H$ is said to be invariant under $T\in \mathcal B(H)$ if $TM\subseteq M$ and reducing if both $M$ and $M^\bot$ are invariant under $T$. For $T\in \mathcal B(H)$, the set $\rho(T):={\{\lambda \in \mathbb C: T-\lambda I:H\rightarrow H\; \text{ is invertible and}\;  (T-\lambda I)^{-1}\in \mathcal B(H) }\}$ is called the resolvent set and the complement $\sigma(T)=\mathbb C\setminus \rho(T)$ is called the spectrum of $T$. The spectral radius of $T$ is given by $m(T)=\sup\{|\lambda|\in \mathbb C: \lambda\in \rho(T)\}.$ It is well known that $\sigma(T)$ is a non empty compact subset of $\mathbb C$. The point spectrum of $T$ is defined by
$\sigma_p(T)={\{\lambda \in \mathbb C: T-\lambda I \; \text{is not one-to-one}}\}.$ Note that $\sigma_{p}(T)\subseteq \sigma(T)$. A self-adjoint operator  $T\in \mathcal B(H)$  is positive if and only if  $\sigma(T)\subseteq [0,\infty)$. If $T\in \mathcal B(H_1,H_2)$, then $T$ is said to be compact if for every  bounded set $S$ of $H_1$, the set $T(S)$ is pre-compact in $H_2$. Similarly, for every bounded sequence $(x_n)$ of $H_1$, $(Tx_n)$ has a convergent subsequence in $H_2$. We denote the set of all compact operators between $H_1$ and $H_2$ by $\mathcal K(H_1,H_2)$. In case if $H_1=H_2=H$, then $\mathcal K(H_1,H_2)$ is denoted by $\mathcal K(H)$. A bounded linear operator $T:H_1\rightarrow H_2$ is called finite rank if $R(T)$ is finite dimensional. The space of all finite rank operators between $H_1$ and $H_2$ is denoted by $\mathcal F(H_1,H_2)$ and we write $\mathcal F(H,H)=\mathcal F(H)$. These standard facts can be obtained in \cite{Hal} and the references therein.

\section{Main Results}
In this section, we give the main results of this work. We begin with the following auxiliary propositions.\\

PROPOSITION 3.1. Let $T\in \mathcal B(H_1,H_2)$ be compact and hyponormal. Then
 \begin{enumerate}
\item [(i).] $m(T)=m(|T|)$
\item[(ii).]  $m(T)=d(0,\sigma(|T|))$
\item [(iii).] $m(T)>0$ if and only if $R(T)$ is closed and $T$ is one-to-one ($T$ is bounded below)
  \item [(iv).] in Particular if $H_1=H_2=H$ and $T^{-1}\in \mathcal B(H)$, then $m(T)=\dfrac{1}{\|T^{-1}\|}$
  \item[(v).] if $H_1=H_2=H$ and $T$ is normal, then
  \begin{itemize}
  \item[(a)] $m(T)=d(0,\sigma(T))$
  \item[(b)]$m(T)=m(T^*)$
  \item[(c)]  $m(T^n)=m(T)^n$ for each $n\in \mathbb N$
  \end{itemize}
 \item [(vi).] if $T\geq 0$, then $m(T)=m(T^{\frac{1}{2}})^2$.
 \end{enumerate}

PROOF. The proof is analogous to the proof of Carvajal and Neves.  See \cite{Car} for proof.

PROPOSITION 3.2. Let $T=K+F+\alpha I$, where $K$ is a positive compact hyponormal operator, $F$ is a self-adjoint finite rank normal operators and $\alpha>0$. Then the following holds:
\begin{enumerate}
 \item [(i).] $R(T)$ is closed
 \item [(ii).] $N(T)$ is finite dimensional. In fact, $N(T)\subseteq R(F)$
 \item [(iii).] $T$ is one-to-one if $K\geq F$
 \item [(iv).] if $T$ is not a finite rank operator, there exists $a>0,\, b>0$ such that $\alpha \in (a\gamma(T), b\gamma(T))$
 \item [(v).] if $T$ is a finite rank operator, then $H$ is finite dimensional.
\end{enumerate}

PROOF. By proposition 3.2 above and  analogous to the proof in \cite{Car} the proof is complete.

PROPOSITION 3.3. Let $T\in \mathcal{B(H)}$  be compact and hyponormal and $\beta\in W_{e}(S)$ where $\alpha >0.$ Then there exists an operator $S\in \mathcal{B(H) }$ such that $\|S\|=\|Z\|,$  $\|S-Z\|<\alpha$ and $T$ is absolutely norm attaining. Furthermore, there exists a vector $\eta\in H,\;\|\eta\|=1$ such that $\|Z\eta\|=\|Z\|$ with $\langle Z\eta,\eta\rangle =\beta.$

PROOF. Consider $S\in \mathcal{B(H) }$ to be contractive then we may assume that $\|S\|=1$ by ignoring the strict inequality. and also that $0<\alpha < 2.$ Let $x_{n}\in H\;(n=1,2,...)$ be such that $\|x_{n}\|=1,$  $\|Sx_{n}\|\rightarrow 1$ and also $\lim_{n\rightarrow\infty}\langle Sx_{n},x_{n}\rangle =\beta.$  Let $S=GL$ be the polar decomposition of $S.$ Here $G$ is a partial isometry and we write $L=\int_{0}^{1}\beta dE_{\beta},$ the spectral decomposition of $L=(S^{*}S)^{\frac{1}{2}}.$ Since $\lim_{n\rightarrow\infty}\|Sx_{n}\|=\|S\|=\|L\|=1,$ we have that $\|Lx_{n}\|\rightarrow 1$ as $n$ tends to $\infty$ and $\lim_{n\rightarrow\infty}\langle Sx_{n},x_{n}\rangle =\lim_{n\rightarrow\infty}\langle GLx_{n},x_{n}\rangle =\lim_{n\rightarrow\infty}\langle Lx_{n},G^{*}x_{n}\rangle.$ Now for $H=\overline{R(L)}\oplus KerL,$ we can choose $x_{n}$ such that $x_{n}\in \overline{R(L)}$ for large $n.$ Indeed, let $x_{n}=x_{n}^{(1)}\oplus x_{n}^{(2)},\; n=1,2,...$ Then we have that $Lx_{n}=Lx_{n}^{(1)}\oplus Lx_{n}^{(2)}=Lx_{n}^{(1)}$ and that $\lim_{n\rightarrow\infty}\|x_{n}^{(1)}\|=1,\;\lim_{n\rightarrow\infty}\|x_{n}^{(2)}\|=0$ since $\lim_{n\rightarrow\infty}\|Lx_{n}\|=1.$ Replacing $x_{n}$ with $\frac{x_{n}^{(1)}}{\|x_{n}^{(1)}\|},$ we obtain
$\lim_{n\rightarrow\infty}\left\|L\frac{1}{\|x_{n}^{(1)}\|}x_{n}^{(1)}\right\|=
\lim_{n\rightarrow\infty}\left\|S\frac{1}{\|x_{n}^{(1)}\|}x_{n}^{(1)}\right\|=1,$
$\lim_{n\rightarrow\infty}\left\langle S\frac{1}{\|x_{n}^{(1)}\|}x_{n}^{(1)},
\frac{1}{\|x_{n}^{(1)}\|}x_{n}^{(1)}\right\rangle =\beta
.$ Now assume that $x_{n}\in \overline{RL}.$ Since $G$ is a partial isometry from $ \overline{R(L)}$ onto $ \overline{R(S)},$ we have that $\|Gx_{n}\|=1$ and $\lim_{n\rightarrow\infty}\langle Lx_{n},G^{*}x_{n}\rangle=\beta.$ Since $L$ is a positive operator, $\|L\|=1$ and for any $x\in H,$ $\langle Lx,x\rangle \leq \langle x,x\rangle =\|x\|^{2}.$ Replacing $x$ with $L^{\frac{1}{2}}x$, we get that $\langle L^{2}x,x\rangle \leq \langle Lx,x\rangle,$ where $L^{\frac{1}{2}}$ is the positive square root of $L.$ Therefore we have that $\|Lx\|^{2}=\langle Lx,Lx\rangle \leq \langle Lx,x\rangle.$ It is obvious that $\lim_{n\rightarrow\infty}\|Lx_{n}\|=1$ and that $\|Lx_{n}\|^{2}\leq\langle Lx_{n},x_{n}\rangle \leq \|Lx_{n}\|^{2}=1.$ Hence,
$\lim_{n\rightarrow\infty}\langle Lx_{n},x_{n}\rangle =1=\|L\|.$ Moreover, Since $I-L\geq 0,$ we have $\lim_{n\rightarrow\infty}\langle (I-L)x_{n},x_{n}\rangle =0.$ thus $\lim_{n\rightarrow\infty}\| (I-L)^{\frac{1}{2}}x_{n}\| =0.$ \\Indeed, $\lim_{n\rightarrow\infty}\| (I-L)x_{n}\| \leq\lim_{n\rightarrow\infty}\| (I-L)^{\frac{1}{2}}\|.\| (I-L)^{\frac{1}{2}}x_{n}\| =0.$ For  $\alpha >0,$ let $\gamma =[0,1-\frac{\alpha}{2}]$ and let $\rho =(1-\frac{\alpha}{2},1].$ We have
$
  L = \int_{\gamma}\mu dE_{\mu}+\int_{\rho}\mu dE_{\mu}
   = LE(\gamma) \oplus LE(\rho).
$
Next we show that $\lim_{n\rightarrow\infty}\|E(\gamma)x_{n}\| =0.$ If there exists a subsequence $x_{n_{i}}, (i=1,2,...,)$ such that $\|E(\gamma)x_{n_{i}}\| \geq\epsilon > 0,\;(i=1,2,...,)$, then since $\lim_{i\rightarrow\infty} \|x_{n_{i}}-Lx_{n_{i}}\| =0,$ it follows from \cite{Oke} that
       $\lim_{i\rightarrow\infty} \|x_{n_{i}}-Lx_{n_{i}}\|^{2} =
       \lim_{i\rightarrow\infty} (\|E(\gamma)x_{n_{i}}-LE(\gamma)x_{n_{i}}\|^{2}
       +\|E(\rho)x_{n_{i}}-LE(\rho)x_{n_{i}}\|^{2})
        = 0.
    $
Hence we have that $\lim_{i\rightarrow\infty} \|E(\gamma)x_{n_{i}}-LE(\gamma)x_{n_{i}}\|^{2} =0.$ Now it is clear that
$ \|E(\gamma)x_{n_{i}}-LE(\gamma)x_{n_{i}}\| \geq \|E(\gamma)x_{n_{i}}\|-\|LE(\gamma)\|.\|E\gamma)x_{n_{i}}\|
  \geq (I-\|LE(\gamma)\|)\|E(\gamma)x_{n_{i}}\|
  \geq \frac{\alpha}{2} \epsilon
  > 0.
$
This is a contradiction. Therefore, $\lim_{n\rightarrow\infty}\|E(\gamma)x_{n}\| =0.$ Since $\lim_{n\rightarrow\infty}\langle Lx_{n},x_{n}\rangle =1,$ we have that
$\lim_{n\rightarrow\infty}\langle LE(\rho)x_{n},E(\rho)x_{n}\rangle =1$ and
$\lim_{n\rightarrow\infty}\langle E(\rho)x_{n},G^{*}E(\rho)x_{n}\rangle =\beta.$
\noindent It is easy to see that $\lim_{n\rightarrow\infty}\| E(\rho)x_{n}\| =1,\;\lim_{n\rightarrow\infty}\left( L\frac{E(\rho)x_{n}}{\|E(\rho)x_{n}\|},\frac{E(\rho)x_{n}}{\|E(\rho)x_{n}\|} \right)=1$ and
$\lim_{n\rightarrow\infty}\left( L\frac{E(\rho)x_{n}}{\|E(\rho)x_{n}\|},G^{*}\frac{E(\rho)x_{n}}{\|E(\rho)x_{n}\|} \right)=\beta$
Replacing $x$ with $\frac{E(\rho)x_{n}}{\|E(\rho)x_{n}\|},$ we can assume that $x_{n}\in E(\rho)H$ for each $n$ and $\|x_{n}\|=1.$ Let
$J = \int_{\gamma}\mu dE_{\mu}+\int_{\rho}\mu dE_{\mu}
   = J_{1} \oplus E(\rho).
$
Then it is evident that $\|J\|=\|S\|=\|L\|=1, Jx_{n}=x_{n}$ and $\|J-L\|\leq\frac{\alpha}{2}.$
If we can find a contraction $V$ such that $V-G\leq\frac{\alpha}{2}$ and $\|Vx_{n}\|=1$ and $\langle Vx_{n},x_{n}\rangle =\beta,$ for a large $n$ then letting $Z=VJ$, we have that $\|Zx_{n}\|=\|VJx_{n}\|=1,$ and that $\langle Zx_{n},x_{n}\rangle =\langle VJx_{n},x_{n}\rangle=\langle Vx_{n},x_{n}\rangle =\beta$
$ \|S-Z\| = \|GL-VJ\| = \alpha.
$ To complete the proof, we now construct the desired contraction $V$.
Clearly, $\lim_{n\rightarrow\infty}\langle x_{n},G^{*}x_{n}\rangle =\beta,$ because
$\lim_{n\rightarrow\infty}\langle L x_{n},G^{*}x_{n}\rangle =\beta$ and
$\lim_{n\rightarrow\infty}\| x_{n}-Lx_{n}\| =0.$ Let $Gx_{n}=\phi_{n}x_{n} +\varphi_{n}y_{n},\;\;(y_{n}\bot x_{n},\;\|y_{n}\|=1)$ then $\lim_{n\rightarrow\infty}\phi_{n} =\beta,$ because $\lim_{n\rightarrow\infty}\langle  Gx_{n},x_{n}\rangle =\lim_{n\rightarrow\infty}\langle  x_{n},G^{*}x_{n}\rangle =\beta$ but $\|Gx_{n}\|^{2}=|\phi_{n}|^{2}+|\varphi_{n}|^{2}=1,$ so we have that
$\lim_{n\rightarrow\infty}|\varphi_{n}| =\sqrt{1-|\beta}|^{2}.$  Then by \cite{Oke} and \cite{Ram} the remaining part of the proof is analogous and this completes the proof.

In this section we consider absolute norm-attainability for commutators of compact hypomormal operators.

LEMMA 3.1. Let $E\in \mathcal{B}(H)$  be compact hyponormal then $EX-XE$  is absolutely norm attaining if there exists a vector $\zeta\in H$ such that $\|\zeta\|=1,\;\;\|E\zeta\|=\|E\|,\;\;\langle E\zeta,\zeta\rangle =0.$

PROOF.  Let $x\in H$ satisfy $x\bot \{\zeta,E\zeta\},$ and define a compact $X$ as follows $X:\zeta \rightarrow \zeta ,\;\;E\zeta \rightarrow -E\zeta, \;\;x\rightarrow 0.$  Since $X$ is a bounded operator on $H$ and $\|X\zeta\|=\|X\|=1,$ $\|EX\zeta -XE\zeta\|=\|E\zeta -(-E\zeta)\|=2\|E\zeta\|=2\|E\|.$ It follows that
$\|EX-XE\|=2\|E\|$ by Proposition 3.1, because $\langle E\zeta,\zeta\rangle =0\in W_{e}(E).$ Hence we have that $\|EX-XE\|=2\|E\|.$ Therefore, $EX-XE$ is absolutely norm attaining.

LEMMA 3.2. Let $S,T\in \mathcal{B}(H)$  be compact hyponormal. If there exists vectors $\zeta, \eta\in H$ such that $\|\zeta\|=\| \eta\|=1,\;\;\|S\zeta\|=\|S\|,\;\|T\eta\|=\|T\|$ and $\frac{1}{\|S\|}\langle S\zeta,\zeta\rangle =- \frac{1}{\|T\|}\langle T\eta,\eta\rangle,$ then $SX-XT$ is is absolutely norm attaining.

PROOF. Since $H$ has an orthonormal basis then by linear dependence of vectors, if $\eta$ and $T\eta$ are linearly dependent, i.e.,$T\eta =\phi\|T\|\eta,$ then we have $|\phi|=1$ and $|\langle T\eta, \eta\rangle |=\|T\|$. It follows that $|\langle S\zeta, \zeta\rangle |=\|S\|$ which implies that $S\zeta =\varphi\|S\|\zeta$ and $|\varphi|=1.$ Hence $\left\langle\frac{S\zeta}{\|S\|},\zeta\right\rangle = \varphi =- \left\langle\frac{T\eta}{\|T\|},\eta\right\rangle =-\phi .$ Defining $X$ as $X:\eta \rightarrow\zeta  ,\;\; \{\eta\}^{\bot} \rightarrow 0,$ we have $\|X\|=1$ and
$(SX-XT)\eta =\varphi(\|S\|+\|T\|)\zeta ,$ which implies that $\|SX-XT\|=\|(SX-XT)\eta\|=\|S\|+\|T\|.$ By \cite{Enf}, it follows that
$\|SX-XT\|=\|S\|+\|T\|=\|\delta_{S,T}\|.$ That is $SX-XT$ is absolutely norm attaining.
\noindent If $\eta$ and $T\eta$ are linearly independent, then $\left|\left\langle\frac{T\eta}{\|T\|},\eta\right\rangle \right|<1,$ which implies that $\left|\left\langle\frac{S\zeta}{\|S\|},\zeta\right\rangle\right|<1.$ Hence $\zeta$ and $S\zeta$ are also  linearly independent. Let us redefine $X$ as follows:
$X:\eta \rightarrow\zeta  ,\;\;\frac{T\eta}{\|T\|}\rightarrow -\frac{S\zeta}{\|S\|},\;\; x \rightarrow 0,$ where $x\in \{\eta,T\eta\}^{\bot}.$
We show that $X$ is a partial isometry. Let $\frac{T\eta}{\|T\|}=\left\langle\frac{T\eta}{\|T\|},\eta\right\rangle \eta +\tau h,\;\; \|h\|=1,\;h\bot\eta .$
Since $\eta$ and $T\eta$ are linearly independent, $\tau\neq 0.$ So we have that
$X\frac{T\eta}{\|T\|}=\left\langle\frac{T\eta}{\|T\|},\eta\right\rangle X\eta +\tau X h=-\left\langle\frac{S\zeta}{\|S\|},\zeta\right\rangle \zeta +\tau X h,$ which implies that
$\left\langle X\frac{T\eta}{\|T\|},\zeta\right\rangle=-\left\langle\frac{S\zeta}{\|S\|},\zeta\right\rangle +\tau \langle X h,\zeta\rangle =-\left\langle\frac{S\zeta}{\|S\|},\zeta\right\rangle . $ It follows then that $\langle X h,\zeta\rangle =0$ i.e., $X h\bot\zeta (\zeta =X\eta).$ Hence we have that $\left\|\left\langle\frac{S\zeta}{\|S\|},\zeta\right\rangle \zeta\right\|^{2}+\|\tau X h\|^{2}=\left\|X\frac{T\eta}{\|T\|}\right\|^{2}
=\left|\left\langle\frac{T\eta}{\|T\|},\eta\right\rangle\right|^{2}+|\tau|^{2}=1,$
which implies that $\|Xh\|=1.$ Now it is evident that $X$ a partial isometry and
$\|(SX-XT)\zeta\|=\|SX-XT\|=\|S\|+\|T\|,$ which is equivalent to $\|\delta_{S,T}(X)\|=\|S\|+\|T\|.$ By Lemma 3.1 and \cite {Oke}, $\|SX-XT\|=\|S\|+\|T\|.$ Hence $SX-XT$ is  absolutely norm attaining.

COROLLARY 3.1.  Let $S,T\in \mathcal{B}(H)$ If both $S$ and $T$ are absolutely norm attaining then the operator  $SXT $ is also absolutely norm attaining.

PROOF.We can assume that $\|S\|=\|T\|=1.$ If both $S$ and $T$ are absolutely norm attaining, then there exists unit vectors $\zeta$ and $\eta$ with $\|S\zeta\|=\|T\eta\|=1.$ We can therefore define
an operator $X$ by $X=\langle \cdot ,T\eta\rangle \zeta$. Clearly, $\|X\|=1.$
Therefore, we have $\|SXT\|\geq \|SXT\eta\|=\|\|T\eta\|^{2}S\zeta\|=1.$
Hence, $\|SXT\|=1,$ that is $SXT$ is also absolutely norm attaining.

PROPOSITION 3.4.  Let $T\in \mathcal {AN}(H)$ be a self-adjoint compact hyponormal operator. Then there exists an orthonormal basis consisting of eigenvectors of $T$.

PROOF. The proof follows in the analogously as in\cite{Car} but we include it for  completeness.
Let $\mathcal B=\{x_\alpha: \alpha\in I\}$ be the maximal set of orthonormal eigenvectors of $T.$ This set is non empty, as $T=T^*\in \mathcal{AN}(H)$. Let $M=\overline{\mbox{span}}\{x_\alpha: \alpha \in I\}$. Then we claim  that $M=H$.
If not, $M^{\bot}$ is a proper non-zero closed subspace of $H$ and it is invariant under $T$.  Since $T=T^* \in \mathcal{AN}(H)$,  then we have either $||T|M^{\bot}||$ or $-||T|M^{\bot}||$ is an eigenvalue for $T|M^{\bot}$. Hence there is a
non-zero vector, say $x_0$ in $M^{\bot}$, such that $Tx_0=\pm ||T|M^{\bot}|| x_0.$ Since $M\cap M^{\bot}={\{0}\},$ we have arrived to a contradiction to the  maximality of $\mathcal B$.

Next, we need to  do a characterization for self-adjoint hyponormal compact operators. We  ask the following question: For a compact hyponormal self-adjoint operator, can we  find $\alpha \in \mathbb R$ such that $K+\alpha I \in \mathcal{AN}(H)$. To solve this first we need
to answer the question when $K+\alpha I\in \mathcal N(H)$. Here we have the following characterization.\\

 LEMMA 3.3. Let $K\in \mathcal K(H)$ be self-adjoint and $a\in \mathbb R$. Let $K=\text{diag}(\lambda_1,\lambda_2,\lambda_3,\dots,)$ with respect to orthonormal basis of $H$. Then the following are equivalent:
 \begin{itemize}
  \item[(i).] $T\in \mathcal N(H)$
  \item[(ii).] there exists $n_0\in \mathbb N$ such that $|\lambda_{n_0}+a|> |a|$.
 \end{itemize}

PROOF. The proof is trivial.

Consider $T=T^*\in \mathcal B(H)$ and have the polar decomposition $T=V|T|$. Let $H_0=N(T),\; H_{+}=N(I-V)$ and $H_{-}=N(I+V)$. Then $H=H_0\oplus H_{+}\oplus H_{-}$ which  are all invariant under $T$. Let $T_0=T|_{N(T)},\; T_{+}=T|_{H_{+}}$ and $T_ {-}=T|_{H_{-}}$. Then $T=T_0\oplus T_{+}\oplus T_{-}$. Further more, $T+$ is strictly positive, $T_{-}$ is strictly negative and $T_0=0$ if $N(T)\neq {\{0}\}$.  Let $P_0=P_{N(T)},\; P_{\pm}=P_{H_{\pm}}$. Then $P_0=I-V^2$ and $P_{\pm}=\frac{1}{2}(V^2\pm V)$. Thus $V=P_{+}-P_{-}$. For details see \cite{Hal}.  

THEOREM 3.1. Let $T\in \mathcal {AN}(H)$ be compact hyponrmal and self-adjoint with the polar decomposition $T=V|T|$. Then
the operator $T$ can be represented as $T=K-F+\alpha V,$
where $K\in \mathcal K(H),\; F\in \mathcal F(H)$ are self-adjoint  with $KF=0$ and $F^2\leq \alpha^2I$

PROOF. Let $H=H_{+}\oplus H_{-}$ and $T=T_{+}\oplus T_{-}$. Since $H_{\pm}$ reduces $T$, we have $T_{\pm}\in \mathcal B(H_{\pm})$. As $T\in \mathcal {AN}(H)$, we have that $T_{\pm}\in \mathcal{AN}(H_{\pm})$.
Hence by  \cite{Enf}, we have that $T_{+}=K_{+}-F_{+}+\alpha I_{H_{+}}$ such that $K_{+}$ is positive compact operator, $F_{+}$ is finite rank positive operator with the property that $K_{+}F_{+}=0$ and $F_{+}\leq \alpha I_{H_{+}} $. As $T_{+}$ is strictly positive, $\alpha>0$. Similarly, $T_{-}\in \mathcal {AN}(H_{-})$ and strictly negative. Hence there exists a triple $(K_{-},F_{-},\beta)$ such that
$-T_{-}=K_{-}-F_{-}+\beta I_{H_{-}}$, where $K_{-}\in \mathcal K(H_{-})$ is positive, $F_{-}\in \mathcal F(H_{-})$ is positive with $K_{-}F_{-}=0,\; F_{-}\leq \beta I_{H_{-}}$ and  $\beta>0$.  The rest follows from \cite{Car} and the proof is complete.

  THEOREM 3.2.  A compact self adjoint hyponormal operator $T\in \mathcal{AN}(H)$ has a countable spectrum.
  
 PROOF. Since $T=T_{+}\oplus T_{-}\oplus T_0$ and all these operators $T_{+},T_{-}$ and $T_0$ are $\mathcal{AN}$ operators. We know that $\sigma(T_{+}),\sigma(T_0)$ are countable, as they are positive. Also, $-T_{-}$ is positive $\mathcal {AN}$-operator and hence $\sigma(T_{-})$ is countable. Hence we can conclude that $\sigma(T)=\sigma(T_{+})\cup \sigma(T_{-})\cup \sigma(T_0)$ is countable.

Now we consider  the structure of normal $\mathcal{AN}$-operators.We see this in the next lemma

LEMMA 3.4. 
 Let $T\in \mathcal{AN}(H)$ be compact hyponormal  with the polar decomposition $T=V|T|$. Then there exists a  compact hyponormal  operator $K$, a  finite rank normal  operator $F\in \mathcal B(H)$ such that $V,K,F$ are mutually commutative.

PROOF. We have $VK=VVK_1=VK_1V=KV$ and $VF=VVF_1=VF_1V=FV$. Also, $KF=0=FK$.

THEOREM 3.3. Let $T\in \mathcal B(H)$ be compact hyponormal.  Then $T\in \mathcal{AN}(H)$ if and only if $T^*\in \mathcal{AN}(H)$.

 PROOF. We know that $T\in \mathcal{AN}(H)$ if and only if $T^*T\in \mathcal{AN}(H)$. Since $T^*T=TT^*$, by Corollary Lemma 3.4 again, it follows that $TT^*\in \mathcal{AN}(H)$ if and only if $T^*\in \mathcal{AN}(H)$.

{\bf Acknowledgement.} This work was partially supported financially by the DFG Grant No. 1603991000.

\end{document}